# AN ANNIHILATING–BRANCHING PARTICLE MODEL FOR THE HEAT EQUATION WITH AVERAGE TEMPERATURE ZERO

By Krzysztof Burdzy[1] and Jeremy Quastel[2]

*University of Washington and University of Toronto*

We consider two species of particles performing random walks in a domain in $\mathbb{R}^d$ with reflecting boundary conditions, which annihilate on contact. In addition, there is a conservation law so that the total number of particles of each type is preserved: When the two particles of different species annihilate each other, particles of each species, chosen at random, give birth. We assume initially equal numbers of each species and show that the system has a diffusive scaling limit in which the densities of the two species are well approximated by the positive and negative parts of the solution of the heat equation normalized to have constant $L^1$ norm. In particular, the higher Neumann eigenfunctions appear as asymptotically stable states at the diffusive time scale.

**1. Introduction.** A branching particle system representation for the heat equation solution with positive temperature was introduced in [4] and later studied in [5] (see also [15] and [16]). Here is an informal description of that model and one of the main results, proved in [5]. Suppose that $D$ is an open set in $\mathbb{R}^d$ and $N$ Brownian particles move independently inside $D$. Whenever one of these particles hits the boundary of $D$, it is killed and one of the other particles, randomly chosen, splits into two particles, so that the number $N$ of particles remains constant. When the number of particles goes to infinity and the initial (normalized) distribution of particles converges to a measure on $D$, then the particle density converges to the solution of the heat equation for every time $t \geq 0$, with the appropriate initial condition, normalized so that it has a constant total mass for all times $t \geq 0$.

The main purpose of this article is to study a related model that involves two different types of particles. We call them + and − particles. The two

Received May 2005; revised November 2005.
[1]Supported in part by NSF Grant DMS-03-03310.
[2]Supported in part by NSERC.

types of particles perform independent symmetric random walks with reflection on the boundary, and annihilate each other on contact. When two particles of different signs annihilate each other, then two other particles are chosen randomly, one + particle and one − particle, and each one of these particles splits into two particles of the same type as the parent particle, so that the numbers of + and − particles are conserved.

The reader might have noticed that we have changed the Brownian particles to symmetric random walks. This is because, except in dimension 1, the Brownian particles would not meet, so one would have to have them annihilate when they came within some $\varepsilon > 0$ of each other. So we might as well just put the particles on a discrete lattice. Hence, particles in our new model will be represented by random walks on the lattice $\varepsilon\mathbb{Z}^d \cap D$, with $\varepsilon^{-d}$ comparable to the number of particles. We will assume that particles reflect from the boundary of $D$ in a way to be described precisely later.

One can easily generalize the model to have more than two types of particles. Particles of different types annihilate on collision and two other particles, from those two types, are chosen at random to split into new particles. The main effect in these systems is segregation of types. At a macroscopic scale one sees the region $D$ decomposed into disjoint regions $D_1(t), \dots, D_m(t)$, each occupied by one of the types. The density of type $i$ evolves on $D_i$ by the heat equation with creation, with Dirichlet boundary conditions on boundaries internal to $D$, and Neumann boundary conditions on boundaries coincident with the boundary of $D$. The regions $D_i$ themselves evolve according to the behavior of the type densities on either side of the boundary. In general, the boundary dynamics are not simple. At the present time the probabilistic picture is based purely on simulations and informal calculations. See [4, 9, 10] and [11] for details. One expects deterministic evolution on the macroscale, as described above. Note that related limiting equations for segregation can be derived in the scaling limit of appropriate systems of reaction–diffusion equations (see [6, 7] and references therein). These problems are quite difficult because, unlike most problems in hydrodynamic limits, the invariant measures of these systems are either unknown or not well understood. Of course, well inside the regions where each individual species resides, one expects a product of Poisson distributions. However, it is not a priori clear that such regions exist, and anyway, all the interesting behavior in the system is at the boundaries between the species.

In this article we consider the special case of two types, and assume, in addition, that the number of each type is the same, say, $N$, where $N$ is proportional to the number of lattice sites, so that the densities are of order 1. Our main result is that, as $N \to \infty$, the densities of the two types of particles converge to the positive and negative parts of the solution $\rho$ of the heat equation normalized to have the correct total masses. Letting $D_+(t)$ and $D_-(t)$ denote the supports of positive and negative parts of the solution, we can see that this corresponds to the general picture just described, with a particular, nontrivial, evolution of the boundary between $D_+$ and $D_-$. The evolution is in the normal direction $\nu$ to the boundary and proportional to $(\Delta + V)\rho_+/\partial_\nu \rho_+ = (\Delta + V)\rho_-/\partial_\nu \rho_-$ there. Here $V = \int_D \delta(\rho = 0)|\nabla\rho|^2\,dx / \int_D |\rho|\,dx$ is the creation term at time $t$ needed to keep the masses conserved.

The method we use is not standard in hydrodynamic limits. The main methods known in hydrodynamic limits are based on entropy or attractiveness. The model is not attractive, and in problems like this one, where the object of interest is a boundary separating two regions on which live mutually singular local equilibrium measures, entropy methods do not seem to be useful. So it is important to develop new methods that can be used for such problems.

The proof is based on an observation that the microscopic evolution of the density $\eta$ of the difference between the occupation numbers of the two types of particles takes a particularly simple form $d\eta = (\Delta + V)\eta\,dt + dM$, where $M$ are martingales, $\Delta$ is the discrete Laplace operator, and $V = V(t)$ is the rate of annihilation. Some calculus and the particular form of the equation allow us to control $V$. Because the proof depends on the particular form of the equation, it does not generalize in an obvious way to systems with more than two types or with unequal numbers of particles. However, it does suggest the precise form the result should take in other, more complicated, cases.

## 2. Preliminaries.

We collect in this section a few results that may be known but we could not find a ready reference for them.

We start with a discussion of a single particle model—reflected random walk on a lattice region approximating a region $D \subset \mathbb{R}^d$. Let $D$ be a bounded connected open set in $\mathbb{R}^d$, $d \geq 2$, and let $D_\varepsilon = \varepsilon\mathbb{Z}^d \cap D$ for $\varepsilon > 0$. We want to construct a nearest-neighbor continuous time random walk $W_t^\varepsilon$ on $D_\varepsilon$. Moreover, we want $W_t^\varepsilon$ to converge to reflecting Brownian motion in $D$ as $\varepsilon \to 0$. A natural choice for $W_t^\varepsilon$ would be a process that jumps to one of its neighbors with equal probabilities, after a suitable time delay. In the main part of our project, such a process is somewhat inconvenient to deal with. We will use a different model for a "reflecting random walk" on $D_\varepsilon$, described below.

The process $W_t^\varepsilon$ is a finite state continuous time Markov process so it is fully described by the jump distribution and holding time at each state. We start by describing its jump distribution. Let $\partial D_\varepsilon$ be the set of $x \in D_\varepsilon$ which have fewer than $2d$ neighbors in $D_\varepsilon$. If $x \in D_\varepsilon \setminus \partial D_\varepsilon$, then the random walk $W_t^\varepsilon$ will jump from $x$ to any of its neighbors with equal probabilities.

We will assume from now on that $D$ has an analytic boundary. This is used later to get some easy estimates on the lattice Laplacian of eigenfunctions of the true Laplacian. One expects the results of this article to hold with much weaker assumptions on the boundary, but we have not pursued this here.

If $\varepsilon$ is very small and $x \in \partial D_\varepsilon$ has fewer than $d$ neighbors in $D_\varepsilon$, then this implies that the normal vector to the boundary close to $x$ is almost parallel to one of the axes. It is not hard to see that by removing from $D_\varepsilon$ all $x \in \partial D_\varepsilon$ with fewer than $d$ neighbors in $D_\varepsilon$, we obtain a set $D'_\varepsilon$ with the property that all $x \in \partial D'_\varepsilon$ have at least $d$ neighbors in $D'_\varepsilon$. By abuse of notation, we will refer to this set as $D_\varepsilon$. Note that every point in $\partial D_\varepsilon$ lies at a distance from $\partial D$ not exceeding $2\varepsilon$.

For $x \in \partial D_\varepsilon$, let $x' \in \partial D$ be the closest point to $x$ on $\partial D$ and let $\mathbf{n}(x')$ be the inward unit normal vector at $x'$. Let $\{\mathbf{e}_1(x') \stackrel{\mathrm{df}}{=} \mathbf{n}(x'), \mathbf{e}_2(x'), \dots, \mathbf{e}_d(x')\}$ be an orthonormal basis, depending on $x$. Let $p_{xy}$ denote the probability that $W^\varepsilon$ makes the next jump from $x$ to $y$ (assuming it is at $x$ now) and note that, by assumption,

$$\sum_{y \in D_\varepsilon, |x-y|=\varepsilon} p_{xy} = 1. \tag{2.1}$$

We want to find $p_{xy}$ for $x \in \partial D_\varepsilon$ so that, for some $c_1 > 0$,

$$\sum_{y \in D_\varepsilon, |x-y|=\varepsilon} p_{xy}(y - x) = c_1 \mathbf{n}(x'). \tag{2.2}$$

Note that equations (2.2) impose only $d-1$ constraints because $c_1 > 0$ is arbitrary. Hence, (2.1) and (2.2) effectively form a set of only $d$ equations, and we have at least $d$ unknowns $p_{xy}$. It is not hard to see that (2.1)–(2.2) have a solution.

Next we choose the holding times. For $x \in D_\varepsilon \setminus \partial D_\varepsilon$, we let the holding time at $x$ have the mean $h_\varepsilon(x) = \varepsilon^2$. This corresponds to the usual space–time scaling for Brownian motion. For $x \in \partial D_\varepsilon$, the mean $h_\varepsilon(x)$ of the holding time at $x$ is chosen so that

$$\lim_{s \downarrow 0} \frac{1}{s} \sum_{1 \le j \le d} E(((W^\varepsilon_{t+s} - W^\varepsilon_t) \cdot \mathbf{e}_j(x'))^2 | W^\varepsilon_t = x) = 1. \tag{2.3}$$

The discrete Laplacian $\Delta_\varepsilon$, that is, the generator of the process $W^\varepsilon$ is given by

$$\Delta_\varepsilon f(x) = h_\varepsilon^{-1}(x) \sum_{y \in D_\varepsilon, |y-x|=\varepsilon} p_{xy}(f(y) - f(x)). \tag{2.4}$$

The following lemma provides an intuitive basis for the interpretation of our main results. Since the lemma is not used later in the paper, we only observe that it can be derived from Theorem 6.3 of [21] and we do not supply a formal proof.

LEMMA 2.1. *Suppose that* $x_\varepsilon \in D_\varepsilon$ *and* $x_\varepsilon \to x_0 \in \overline{D}$ *as* $\varepsilon \to 0$. *If* $W_0^\varepsilon = x_\varepsilon$, *then* $\{W_t^\varepsilon, t \geq 0\}$ *converge weakly to the reflecting Brownian motion on* $\overline{D}$ *starting from* $x_0$, *as* $\varepsilon \to 0$.

We have assumed that $D$ has an analytic boundary. This implies that $D$ has a discrete spectrum for the Laplacian with Neumann boundary conditions; see, for example, Section 2 of [1]. In the Dirichlet case, the spectrum is discrete for any bounded open connected set $D$. In the Neumann case, the spectrum is not necessarily discrete if $D$ is bounded but has a nonsmooth boundary; see, for example, [17]. It is quite easy to see that the spectrum is discrete for the Neumann Laplacian when $\partial D$ can be represented locally by a Lipschitz function; see, for example, Section 2 of [1]. This condition on the boundary can be substantially relaxed—it is enough to suppose that the boundary is locally represented by a function which is Hölder continuous with sufficiently small exponent, see, for example, [19]. As we said, our assumptions on the smoothness of $\partial D$ are much stronger than that, for technical reasons.

Let $-\lambda_n$, $n = 0, 1, 2, \ldots$, be the eigenvalues of the Laplacian on $D$ with the Neumann boundary conditions, and let $\phi_n$ be the corresponding eigenfunctions. We list the same eigenvalue more than once, if it has a multiplicity greater than 1. For a measure $\mu$ on $\overline{D}$, we let $\hat{\mu}_n = \int_{\overline{D}} \phi_n \, d\mu$.

The following review of random measures is based on Section 3 of [12]. A sequence of finite measures $\mu^n$ on a space $\Lambda$ converges weakly to $\mu$ if for any bounded continuous function $f$, we have $\int_\Lambda f(x)\, d\mu^n(x) \to \int_\Lambda f(x)\, d\mu(x)$ as $n \to \infty$. Let $M_F(\Lambda)$ denote the space of all finite measures on $\Lambda$ equipped with the topology of weak convergence, let $M_1(\Lambda)$ be the subspace of $M_F(\Lambda)$ consisting of probability measures, and let $M_{F,c}(\Lambda)$ denote the subspace of $M_F(\Lambda)$ consisting of measures with total variation less than or equal to $c$. Let $\mathcal{B}(\overline{D})$ denote the family of Borel subsets of $\overline{D}$. Suppose that, for some probability space $(\Omega, \mathcal{F}, P)$, the function $\mu : \mathcal{B}(\overline{D}) \times \Omega \to \mathbb{R}$ has the following two properties: (i) for a fixed $\omega$, $\mu(\cdot, \omega)$ is a finite measure on $\overline{D}$, and (ii) for a fixed $A \in \mathcal{B}(\overline{D})$, $\mu(A, \cdot)$ is a random variable. Then the distribution of $\mu$ is an element of $M_1(M_F(\overline{D}))$. We will refer to $\mu$ as a "random measure." If $\Lambda$ is a metric space, then the space of right-continuous functions with left limits, $f : (0, \infty) \to \Lambda$, equipped with the Skorohod topology (see Section 3.6 of [12]) will be denoted by $S((0, \infty), \Lambda)$. We will use the space $M_1(S((0, \infty), M_F(\overline{D})))$ to state our main results. To give a meaning to this symbol, we have to specify a metric on $M_F(\overline{D})$. We will use the Prohorov metric which induces a topology on $M_F(\overline{D})$ equivalent to weak convergence of measures (see [13], Chapter 3, Sections 1 and 3). In this article we will be concerned with the convergence of processes on the open half-line $(0, \infty)$, not the usual semi-closed half-line $[0, \infty)$. The convergence in the Skorohod topology on $(0, \infty)$ is defined as the convergence in the Skorohod topology on every compact subinterval of $(0, \infty)$.

LEMMA 2.2. *Suppose that* $f:\overline{D}\to\mathbb{R}$ *is a continuous (and, hence, bounded) function and let* $a_n=\int_{\overline{D}}f(x)\phi_n(x)\,dx$. *Suppose that* $c_1<\infty$ *and* $\mu^k$, $k\ge 1$, *are random measures (possibly defined on different probability spaces), with distributions in* $M_1(M_{F,c_1}(\overline{D}))$. *Assume that, for every fixed* $n$, $\lim_{k\to\infty}\hat{\mu}_n^k=a_n$ *in distribution. Then the distributions of* $\mu^k$ *converge weakly in* $M_1(M_F(\overline{D}))$ *to* $\delta_\mu$, *where* $\mu(dx)=f(x)\,dx$.

PROOF. Let $d\nu^k(x)=d\mu^k(x)-f(x)\,dx$ and note that $\hat{\nu}_n^k\to 0$ in distribution as $k\to\infty$, for every $n$. It will suffice to show that the distributions of $\nu^k$ converge weakly in $M_1(M_F(\overline{D}))$ to $\delta_{\mathbf{0}}$, where $\mathbf{0}$ is the measure identically equal to 0. Since $\overline{D}$ is compact, so is $M_{F,c_1}(\overline{D})$ and it follows that the sequence $\nu^k$ is tight and contains a convergent subsequence. Let $\nu$ be the weak limit of a subsequence of $\nu^{k_j}$. It will be enough to show that $\nu=\delta_{\mathbf{0}}$.

Let $G_n:M_F(\overline{D})\to\mathbb{R}$ be defined by $G_n(\sigma)=|\int_{\overline{D}}\phi_n(x)\,d\sigma(x)|\wedge 1$. The functionals $G_n$ are continuous and bounded, and distributions of $\nu^{k_j}$ converge weakly to $\nu$ in $M_1(M_F(\overline{D}))$, so for every fixed $n$, $EG_n(\nu^{k_j})\to EG_n(\nu)$ as $j\to\infty$. By assumption, $\hat{\nu}_n^k\to 0$ in distribution for every $n$, so $EG_n(\nu)=\lim_{j\to\infty}EG_n(\nu^{k_j})=0$ for every $n$. Hence, $\nu$ is supported on measures $\sigma\in M_F(\overline{D})$ with the property that $\hat{\sigma}_n=0$ for all $n$. It will suffice to show that every measure with this property and finite total variation is identically equal to 0.

Fix a nonrandom measure $\sigma$ on $\overline{D}$ with a finite total variation and such that $\hat{\sigma}_n=0$ for all $n$. Fix any Borel set $A\subset\overline{D}$. It will be enough to show that $\int_{\overline{D}}\mathbf{1}_A(x)\,d\sigma(x)=0$.

According to the Weyl formula, $\lambda_n\sim n^{2/d}$ (see [8], Vol. I, Chapter VI, Section 4.4; or see [19] for a recent strong version of this theorem). By Theorem 1 of [14], $\|\phi_n\|_\infty\le c_2\lambda_n^{(d-1)/4}$, so $\|\phi_n\|_\infty\le c_3 n^{(d-1)/2d}$.

Let $P_t$ be the transition semigroup for the reflected Brownian motion in $\overline{D}$ and for some fixed $t>0$, let $g(x)=P_t\mathbf{1}_A(x)$. Then $g(x)=\sum_n(\hat{\mathbf{1}}_A)_n\phi_n(x)e^{-\lambda_n t}$, where $(\hat{\mathbf{1}}_A)_n=\int_{\overline{D}}\mathbf{1}_A(x)\phi_n(x)\,dx$. If we write $g(x)=\sum_n\hat{g}_n\phi_n(x)$, then

$$|\hat{g}_n|\le|D|\|\phi_n\|_\infty e^{-\lambda_n t}\le c_4 n^{(d-1)/2d}e^{-c_5 n^{2/d}t}.$$

Without loss of generality, assume that the total variation of $\sigma$ is not greater than 1 and note that

$$\int_{\overline{D}}|\phi_n(x)|\,d|\sigma(x)|\le|D|\|\phi_n\|_\infty\le c_6 n^{(d-1)/2d}.$$

This and other estimates obtained so far imply that

$$\sum_n\int_{\overline{D}}|\hat{g}_n||\phi_n(x)|\,d|\sigma(x)|\le\sum_n c_4 n^{(d-1)/2d}e^{-c_5 n^{2/d}t}c_6 n^{(d-1)/2d}<\infty.$$

Hence, we can change the order of integration and summation in the following formula:

$$\int_D g(x)\,d\sigma(x) = \int_D \sum_n \hat{g}_n \phi_n(x)\,d\sigma(x) = \sum_n \hat{g}_n \int_D \phi_n(x)\,d\sigma(x) = 0.$$

We have proved that $\int_D P_t \mathbf{1}_A(x)\,d\sigma(x) = 0$ for every $t > 0$. Clearly, $|P_t \mathbf{1}_A(x)| \leq 1$ for all $x$, and $P_t \mathbf{1}_A(x) \to \mathbf{1}_A(x)$ for almost every $x \in \overline{D}$, so by the dominated convergence,

$$\int_D \mathbf{1}_A(x)\,d\sigma(x) = \lim_{t\to 0} \int_D P_t \mathbf{1}_A(x)\,d\sigma(x) = 0. \qquad \square$$

The following result can be proved using standard methods, so we leave its proof to the reader.

LEMMA 2.3. *Suppose that, for every $\varepsilon > 0$, we have a real-valued process $\{R^\varepsilon(t), t > 0\}$ which is equal to $R_1^\varepsilon(t) + R_2^\varepsilon(t) + R_3^\varepsilon(t)$, and these processes satisfy the following conditions:*

(i) *For every $\varepsilon > 0$, $t \to R_1^\varepsilon(t)$ is right-continuous and nondecreasing a.s.*
(ii) *For every fixed $t > 0$, the family $\{R_1^\varepsilon(t), \varepsilon > 0\}$ is tight.*
(iii) *For every fixed $t > 0$, $\sup_{0\leq s\leq t} |R_2^\varepsilon(s)|$ converges to $0$ in distribution, as $\varepsilon \to 0$.*
(iv) *For every $t > 0$,*

$$\limsup_{\delta_1,\delta_2\to 0} \limsup_{\varepsilon\to 0} P\left( \sup_{0\leq t_1,t_2\leq t, |t_1-t_2|\leq\delta_1} |R_3^\varepsilon(t_2) - R_3^\varepsilon(t_1)| \geq \delta_2 \right) = 0.$$

*Then the family $\{R^\varepsilon, \varepsilon > 0\}$ is tight in $M_1(S((0,\infty),\mathbb{R}))$.*

**3. Main results.** We will now describe the main object of our study, a particle system. Our description will be partly informal because this will make it more accessible without loss of rigor.

The state of the particle system at time $t \geq 0$ will be encoded as an integer-valued random function $\eta(t) = \eta_x(t)$ on $D_\varepsilon$. We will often suppress $t$ in the notation and write $\eta_x$, with the convention that if $\eta_x > 0$, then there are $\eta_x$ particles of type $+$ at $x \in D_\varepsilon$, and $\eta_x < 0$ signifies the presence of $|\eta_x|$ particles of type $-$ at $x$. Obviously, $\eta_x = 0$ means that there are no particles at $x$. Since $+$ and $-$ particles annihilate each other instantaneously, we do not have to have a notation representing both types of particles at the same site of $D_\varepsilon$. We assume that there are $N$ particles of type $+$ and $N$ particles of type $-$ for every $t$.

The easiest way to describe the evolution of $\eta$ is to use the particle picture. Each particle performs continuous time symmetric simple random walk (defined in Section 2), independent of other particles, until one of the particles hits a particle of the other type. When a particle jumps to a site occupied by at least one particle of the different sign, two particles of opposite signs annihilate each other. At the same time, one particle is chosen uniformly from the family of + particles and one − particle is chosen uniformly as well. Each of the two chosen particles splits into two offspring of the same type as the parent, so that the number of particles of each type remains constant. The particles move independently until the next annihilation and birth event, and the evolution continues in the same manner. From the point of view of the mathematical description of the model, it is more convenient to represent the "annihilation and birth" event as a jump of the annihilated + particle to a randomly chosen + particle, and the same for the annihilated − particle.

In the above description of the dynamics of the system, when we say that a particle is chosen "uniformly," it means that we choose one of the $N$ particles with the same probability $1/N$; in other words, a particle that is annihilated may be the one that splits into two offspring. In such a case, the offspring are born at the site where the annihilated particle resided just before the jump that lead to its annihilation. This is a different convention than in [5], but this convention will simplify some formulas and, of course, it makes little difference when $N$ is large.

The informal description given above can be rigorously expressed in terms of the generator $L = L_\varepsilon$ for the process. Let $a^+ = \max(a, 0)$ and $a^- = \max(-a, 0)$. The configuration $\eta$ such that $\eta_x = 1$ and $\eta_y = 0$ for $y \neq x$ will be denoted $\mathbf{I}_x$. The formula for the generator $L$ of the process $\eta$ is

$$
\begin{aligned}
&Lf(\eta)\\
&\quad = \sum_{x,y\in D_\varepsilon} h_\varepsilon^{-1}(x) p_{xy} \Bigg\{ \eta_x^+ \mathbf{1}_{\{\eta_y \geq 0\}} (f(\eta - \mathbf{I}_x + \mathbf{I}_y) - f(\eta))\\
&\qquad\qquad + \eta_x^- \mathbf{1}_{\{\eta_y \leq 0\}} (f(\eta + \mathbf{I}_x - \mathbf{I}_y) - f(\eta))\\
&\qquad\qquad + \eta_x^+ \mathbf{1}_{\{\eta_y < 0\}} N^{-2}\\
&\qquad\qquad\quad \times \sum_{u,v\in D_\varepsilon} \eta_u^+ \eta_v^- (f(\eta - \mathbf{I}_x + \mathbf{I}_y + \mathbf{I}_u - \mathbf{I}_v) - f(\eta))\\
&\qquad\qquad + \eta_x^- \mathbf{1}_{\{\eta_y > 0\}} N^{-2}\\
&\qquad\qquad\quad \times \sum_{u,v\in D_\varepsilon} \eta_u^- \eta_v^+ (f(\eta + \mathbf{I}_x - \mathbf{I}_y - \mathbf{I}_u + \mathbf{I}_v) - f(\eta)) \Bigg\}.
\end{aligned}
$$

We would like to point out that the first sum in the above formula extends over only those $x, y \in D_\varepsilon$ that satisfy $|x - y| = \varepsilon$. This extra condition is enforced automatically by the presence of the factor $p_{xy}$. A similar remark applies to other formulas in this paper.

The normalized particle density $u(x,t) = u^{N,\varepsilon}(x,t)$ is defined by $u(x,t) = N^{-1}\varepsilon^{-d}\eta_x$ for $x \in D_\varepsilon$. Typically, we will be interested in the population size of order $N \sim \varepsilon^d$. Let $P_t$ denote the heat semigroup on $\overline{D}$, corresponding to the reflecting Brownian motion, and for a measure $\mu$ on $\overline{D}$, let $P_t\mu$ denote the measure with the density $\int_{\overline{D}} P_t(x,\cdot)\,d\mu(x)$. When $\mu = \sum_{x\in D_\varepsilon} \varepsilon^d u^{N,\varepsilon}(x,0)\mathbf{i}_x$ and $\mathbf{i}_x$ denotes the probability measure with the unit mass at $x$ (by abuse of notation), then we will write $P_t u^{N,\varepsilon}(0)$ to denote $P_t\mu$. In other words, $P_t u^{N,\varepsilon}(0)\,dy = \sum_{x\in D_\varepsilon} u^{N,\varepsilon}(x,0)P_t(x,y)\,dy$.

For a nonzero measure $\mu$ on $\overline{D}$ of finite total variation, we define $\overline{\mu}$ to be $c\mu$, where $c = c(\mu)$ is a normalizing constant chosen so that the total variation of $\overline{\mu}$ is equal to 2. For definiteness, we let $\overline{\mu} \equiv 0$ if $\mu \equiv 0$.

For the meaning of $M_1(S((0,\infty), M_F(\overline{D})))$, see Section 2.

THEOREM 3.1. *Suppose that $\varepsilon \to 0$ and $N \to \infty$ in such a way that $c_1\varepsilon^{-d} \le N \le c_2\varepsilon^{-d}$ for some constants $0 < c_1 < c_2 < \infty$. Assume that $D$ has an analytic boundary, $u^{N,\varepsilon}(x,0)$ are nonrandom, and for some signed measure $\mu$ on $\overline{D}$ which does not vanish identically, $\sum_{x\in D_\varepsilon} \varepsilon^d u^{N,\varepsilon}(x,0)\mathbf{i}_x \to \mu$ in $M_F(\overline{D})$ as $\varepsilon \to 0$. Then, as $\varepsilon \to 0$,*

$$\sum_{x\in D_\varepsilon} \varepsilon^d u^{N,\varepsilon}(x,\cdot)\mathbf{i}_x - \overline{P_\cdot\mu} \to \delta_{\mathbf{0}} \qquad \text{in } M_1(S((0,\infty), M_F(\overline{D}))).$$

The symbol $\delta_{\mathbf{0}}$ in the last formula refers to the process identically equal to the null measure.

Note that we have not assumed that the total variation of $\mu$ is equal to 2. The total variation of $\mu$ cannot be larger than 2, but it can be smaller than 2 when the particles are tightly interspersed at time 0. If the total variation of $\mu$ is less than 2, then for small $\varepsilon$, the particle configuration has an almost instantaneous jump at time 0 to a configuration that approximates $\overline{\mu}$. For this reason, we obtain convergence only in $M_1(S((0,\infty), M_F(\overline{D})))$, not in $M_1(S([0,\infty), M_F(\overline{D})))$. If the total variation of $\mu$ is equal to 2, our arguments show that the convergence holds in $M_1(S([0,\infty), M_F(\overline{D})))$.

The assumption that $u^{N,\varepsilon}(x,0)$ are nonrandom measures is easy to remove—we added it for technical convenience only.

Theorem 3.1 is a special case of Theorem 3.2 below. We need some more notation to present this more general result.

Recall from Section 2 that $-\lambda_n$, $n = 0,1,2,\dots$, are the eigenvalues of the Laplacian on $D$ with the Neumann boundary conditions, and $\phi_n$ are the corresponding eigenfunctions. For a measure $\mu$ on $\overline{D}$, we write $\hat{\mu}_n = \int_{\overline{D}} \phi_n\,d\mu$. For a function $f : D_\varepsilon \to \mathbb{R}$, we let $\langle f, g\rangle = \varepsilon^d \sum_{x\in D_\varepsilon} f(x)g(x)$. The Fourier coefficients for the "density" $u(x,t)$ will be denoted $\hat{u}_n = \hat{u}_n(t) = \langle u(t), \phi_n\rangle$. In other words, $\hat{u}_n$ is the $n$th Fourier coefficient for the measure $\sum_{x\in D_\varepsilon} \varepsilon^d u^{N,\varepsilon}(x,t)\mathbf{i}_x$.

Note that if $\sum_{x\in D_\varepsilon} \varepsilon^d u^{N,\varepsilon}(x,0)\mathbf{i}_x \to \mu$ in $M_F(\overline{D})$ as $\varepsilon \to 0$, where $\mu$ is a signed measure that is not identically equal to 0, then, for some $n$, there exists $a > 0$ such that $|\hat{u}_n(0)| \geq a$ for sufficiently small $\varepsilon > 0$ (see Lemma 2.2 and its proof).

THEOREM 3.2. *Suppose that $\varepsilon \to 0$ and $N \to \infty$ in such a way that $c_1\varepsilon^{-d} \leq N \leq c_2\varepsilon^{-d}$ for some constants $0 < c_1 < c_2 < \infty$. Assume that $D$ has an analytic boundary and for some $n$, there exists $a > 0$ such that $\inf_{N,\varepsilon} |\hat{u}_n(0)| = a$. Then, as $\varepsilon \to 0$,*

$$\sum_{x\in D_\varepsilon} \varepsilon^d u^{N,\varepsilon}(x,\cdot)\mathbf{i}_x - \overline{P.u^{N,\varepsilon}(0)} \to \delta_{\mathbf{0}} \qquad \text{in } M_1(S((0,\infty), M_F(\overline{D}))).$$

The remainder of the paper is devoted to the proof of this theorem. We start with a very informal overview. The proof will be divided into several steps.

In Step 1 we observe that the density field, $u$, of + particles minus − particles, satisfies a particularly simple, linear equation at the microscopic scale,

$$du_\epsilon = [\Delta_\epsilon^* + V_\epsilon]u_\epsilon + dM_\epsilon,$$

where $\Delta_\epsilon^*$ is the adjoint of the random walk generator, $V = V(t)$ is the particle annihilation rate, and $M$ is a field of martingales.

This is based on the following elementary observation. Let $A$ be the generator of a continuous time Markov process with state space $S$. Let $L$ denote the generator of a system of particles of two types, on $S$, performing this dynamics, and, in addition, annihilating on contact. Let $\eta_x$ denote the number of the first type minus the number of the second type, at $x \in S$. Then

$$L(\eta_x) = (A^*\eta)_x.$$

Adding the particle creation term, which is clearly of mean field type, gives the preceding linear equation.

Now our problem is to show that $u_\epsilon(t,x)$ is close to $v(t,x) = \overline{P_t u_\epsilon(0,x)}$. $v(t,x)$ satisfies a similar looking equation $\partial_t v = [\Delta + W]v$, where the job of $W = W(t)$ is to maintain the total mass of $v$. Formally differentiating $\int_D |v|\,dx$ gives $W(t) = \int_D \delta(v=0)|\nabla v|^2\,dx/\int_D |v|\,dx$. So there are two key things to prove: First, that the martingale terms vanish in the limit. Second, that $V_\epsilon$ looks like $W$ on the macroscale. Note that control of $M_\epsilon$ itself is not enough. One really needs to control martingale terms like $\int_0^t e^{\int_s^t V_\epsilon(u)\,du}\,dM_\epsilon(s)$. The $V_\epsilon$ are the main unknowns, and this is the key difficulty of the problem.

The linear form of the equation suggests the use of Fourier analysis. In Step 2 we rewrite the equation as a system of stochastic differential equations for the Fourier coefficients of the density field [see (3.6) and (3.7)]. A simple observation is that, because of the conservation law, and well-known bounds on the eigenfunctions of the Laplacian, the Fourier coefficients themselves are bounded. A bit of calculus, and the particular form of the system, then allow us to obtain an a priori estimate on $V_\epsilon$ [see (3.15)]. Once this is done, a preliminary form of the limiting equation can be obtained [see (3.16)].

(3.16) tells us that what we see macroscopically at time $t$ is the evolution of the density field by the heat equation, as long as we are willing to multiply by some scaling factor. At first glance, this would appear to imply the full result, because the scaling factor is fixed by the conservation law. However, there is still a lot of work to show that this is in fact the case. What could be happening is that the two type of particles are coexisting on some mesoscopic scale. All we would see on the macroscale is a net decrease in the total mass. In Step 3 we show that this cannot happen. The reason is that such a situation would lead to a very large rate of annihilation, as + and − particles would find themselves unnaturally close to each other. So $V_\epsilon$ would get very large. But we already have a bound (3.15) which prevents this. Finally, in Steps 4 and 5 similar ideas are used to show the required tightness.

We now proceed with the proof.

Proof of Theorem 3.2. *Step* 1. In this step we will show that, for $z \in D_\varepsilon$ and $f(\eta) = \eta_z$,

$$L_\varepsilon f(\eta) = \Delta_\varepsilon^* \eta_z + V \eta_z, \tag{3.1}$$

where $V = V(\varepsilon, N, \eta)$ is the (normalized, instantaneous) jump intensity for the particle system in state $\eta$, defined by

$$V = 2N^{-1} \sum_{x,y \in D_\varepsilon} h_\varepsilon^{-1}(x) p_{xy} (\eta_x^+ \mathbf{1}_{\{\eta_y < 0\}} + \eta_x^- \mathbf{1}_{\{\eta_y > 0\}})$$

and

$$\Delta_\varepsilon^* f(x) = \sum_{y \in D_\varepsilon} (h_\varepsilon^{-1}(y) p_{yx} f(y) - h_\varepsilon^{-1}(x) p_{xy} f(x)).$$

This operator is the adjoint of the discrete Laplacian $\Delta_\varepsilon$ given by

$$\Delta_\varepsilon f(x) = h_\varepsilon^{-1}(x) \sum_{y \in D_\varepsilon} p_{xy} (f(y) - f(x)).$$

Fix some $z \in D_\varepsilon$ and let $f(\eta) = \eta_z$. Then

$$\begin{aligned} Lf(\eta) = \sum_{x,y \in D_\varepsilon} h_\varepsilon^{-1}(x) p_{xy} \Bigg\{ & \eta_x^+ \mathbf{1}_{\{\eta_y \geq 0\}} (\mathbf{1}_{\{y=z\}} - \mathbf{1}_{\{x=z\}}) \\ & + \eta_x^- \mathbf{1}_{\{\eta_y \leq 0\}} (\mathbf{1}_{\{x=z\}} - \mathbf{1}_{\{y=z\}}) \\ & + \eta_x^+ \mathbf{1}_{\{\eta_y < 0\}} N^{-2} \end{aligned}$$

$$
\begin{aligned}
&\times \sum_{u,v\in D_\varepsilon} \eta_u^+\eta_v^- (\mathbf{1}_{\{y=z\}} + \mathbf{1}_{\{u=z\}} \\
&\qquad - \mathbf{1}_{\{x=z\}} - \mathbf{1}_{\{v=z\}}) \\
&+ \eta_x^- \mathbf{1}_{\{\eta_y>0\}} N^{-2} \\
&\times \sum_{u,v\in D_\varepsilon} \eta_u^-\eta_v^+ (\mathbf{1}_{\{x=z\}} + \mathbf{1}_{\{v=z\}} \\
&\qquad - \mathbf{1}_{\{y=z\}} - \mathbf{1}_{\{u=z\}}) \Bigg\}.
\end{aligned} \tag{3.2}
$$

The sum of the terms in (3.2) with the indicators $\mathbf{1}_{\{x=z\}}$ is equal to

$$
\begin{aligned}
&\sum_{y\in D_\varepsilon} h_\varepsilon^{-1}(z) p_{zy} \Bigg[ -\eta_z^+ \mathbf{1}_{\{\eta_y\geq 0\}} + \eta_z^- \mathbf{1}_{\{\eta_y\leq 0\}} \\
&\qquad - \eta_z^+ \mathbf{1}_{\{\eta_y<0\}} N^{-2} \sum_{u,v\in D_\varepsilon} \eta_u^+\eta_v^- \\
&\qquad + \eta_z^- \mathbf{1}_{\{\eta_y>0\}} N^{-2} \sum_{u,v\in D_\varepsilon} \eta_u^-\eta_v^+ \Bigg] \\
&= \sum_{y\in D_\varepsilon} h_\varepsilon^{-1}(z) p_{zy} [ -\eta_z^+ \mathbf{1}_{\{\eta_y\geq 0\}} + \eta_z^- \mathbf{1}_{\{\eta_y\leq 0\}} \\
&\qquad - \eta_z^+ \mathbf{1}_{\{\eta_y<0\}} N^{-2} N^2 + \eta_z^- \mathbf{1}_{\{\eta_y>0\}} N^{-2} N^2 ] \\
&= - \sum_{y\in D_\varepsilon} h_\varepsilon^{-1}(z) p_{zy} \eta_z.
\end{aligned} \tag{3.3}
$$

A similar calculation shows that the sum of the terms in (3.2) with the indicators $\mathbf{1}_{\{y=z\}}$ is equal to

$$
\sum_{x\in D_\varepsilon} h_\varepsilon^{-1}(x) p_{xz} \eta_x. \tag{3.4}
$$

The sum of the terms in (3.2) with the indicators $\mathbf{1}_{\{u=z\}}$ is equal to

$$
\begin{aligned}
&\sum_{x,y\in D_\varepsilon} h_\varepsilon^{-1}(x) p_{xy} \Bigg[ \eta_x^+ \mathbf{1}_{\{\eta_y<0\}} N^{-2} \sum_{v\in D_\varepsilon} \eta_z^+\eta_v^- \\
&\qquad - \eta_x^- \mathbf{1}_{\{\eta_y>0\}} N^{-2} \sum_{v\in D_\varepsilon} \eta_z^-\eta_v^+ \Bigg] \\
&= \sum_{x,y\in D_\varepsilon} h_\varepsilon^{-1}(x) p_{xy} [\eta_x^+ \mathbf{1}_{\{\eta_y<0\}} N^{-2} \eta_z^+ N - \eta_x^- \mathbf{1}_{\{\eta_y>0\}} N^{-2} \eta_z^- N]
\end{aligned} \tag{3.5}
$$

$$= (1/2)V\eta_z.$$

Similarly, the sum of the terms in (3.2) with the indicators $\mathbf{1}_{\{v=z\}}$ is equal to $(1/2)V\eta_z$. Combining this and (3.2)–(3.5), we obtain (3.1).

*Step* 2. We will now derive estimates for the Fourier coefficients of $u$. Note that $\eta$ is a finite state continuous time Markov process. Therefore (see [13], Chapter 4, Proposition 1.7 and Chapter 4, Section 2), $M_n(t) \stackrel{\text{df}}{=} \hat{u}_n(t) - \int_0^t L\hat{u}_n(s)\,ds$ is a martingale. In other words, the Fourier coefficient $\hat{u}_n$ satisfies the following stochastic differential equation:

$$d\hat{u}_n = L\hat{u}_n\,dt + dM_n. \tag{3.6}$$

Since $\phi_n$ is the $n$th eigenfunction, $\Delta\phi_n = -\lambda_n\phi_n$, where $\Delta$ is the Laplacian on $D$ with the Neumann boundary conditions. Recall that we have assumed that $D$ has an analytic boundary. It was pointed out in [2] that, in view of Theorem 5.7.1 on page 169 of [18], if the coefficients of a second-order elliptic equation are real analytic on a bounded analytic domain $D$ up to the boundary, then, for every point $z$ on $\partial D$, there exists a ball $B$ centered at $z$ such that solutions of the elliptic equation can be extended to be real analytic functions on $B$. Hence, for every point $x \in \overline{D}$, $\phi_n(x)$ is equal to its power series in a neighborhood of $x$. The radius of convergence of the series is strictly positive for every $x \in \overline{D}$ and continuous as a function of $x$. Since the domain $D$ is bounded, the radius of convergence is bounded below for all $x \in \overline{D}$ by a strictly positive constant. The coefficients of the power series are continuous as functions of $x$. We will estimate $\sum_{y\in D_\varepsilon} p_{xy}(\phi_n(y) - \phi_n(x))$ using power series expansion at $x$. Since the normal derivative of $\phi_n$ vanishes on the boundary, (2.2) implies that the sum of linear terms vanishes. The sum of quadratic terms is equal to $\varepsilon^2\Delta\phi_n(x) + o(\varepsilon^3)$ because of (2.3). The contribution of higher terms is $o(\varepsilon^3)$. Hence, the sum in question is equal to $\varepsilon^2\Delta\phi_n(x) + o(\varepsilon^3)$. It is easy to check using the definition (2.3) that $h_\varepsilon(x)$ is of order $\varepsilon^2$. All these estimates imply that we have

$$\Delta_\varepsilon\phi_n(x) = \Delta\phi_n(x) + \psi_\varepsilon(x) = -\lambda_n\phi_n(x) + \psi_\varepsilon(x),$$

where $|\psi_\varepsilon(x)| \le c(n)\varepsilon$. We have $L\eta_x = \Delta_\varepsilon^*\eta_x + V\eta_x$ and

$$\hat{u}_n = \langle u, \phi_n\rangle = \varepsilon^d \sum_{x\in D_\varepsilon} N^{-1}\varepsilon^{-d}\eta_x\phi_n(x) = N^{-1}\sum_{x\in D_\varepsilon}\eta_x\phi_n(x),$$

so

$$\begin{aligned}
L\hat{u}_n &= N^{-1}\sum_{x\in D_\varepsilon}(L\eta_x)\phi_n(x)\\
&= N^{-1}\sum_{x\in D_\varepsilon}(\Delta_\varepsilon^*\eta_x)\phi_n(x) + N^{-1}\sum_{x\in D_\varepsilon}(V\eta_x)\phi_n(x)
\end{aligned}$$

$$= N^{-1} \sum_{x\in D_\varepsilon} \left( \sum_{y\in D_\varepsilon} (h_\varepsilon^{-1}(y) p_{yx} \eta_y - h_\varepsilon^{-1}(x) p_{xy} \eta_x) \right) \phi_n(x)$$

$$+ N^{-1} V \sum_{x\in D_\varepsilon} \eta_x \phi_n(x)$$

$$= N^{-1} \sum_{x\in D_\varepsilon} \left( \sum_{y\in D_\varepsilon} (h_\varepsilon^{-1}(x) p_{xy} \phi_n(y) - h_\varepsilon^{-1}(x) p_{xy} \phi_n(x)) \right) \eta_x \tag{3.7}$$

$$+ N^{-1} V \sum_{x\in D_\varepsilon} \eta_x \phi_n(x)$$

$$= N^{-1} \sum_{x\in D_\varepsilon} \eta_x (\Delta_\varepsilon \phi_n(x)) + V N^{-1} \sum_{x\in D_\varepsilon} \eta_x \phi_n(x)$$

$$= N^{-1} \sum_{x\in D_\varepsilon} \eta_x (-\lambda_n \phi_n(x) + \psi_\varepsilon(x)) + V N^{-1} \sum_{x\in D_\varepsilon} \eta_x \phi_n(x)$$

$$= (V - \lambda_n) \hat{u}_n + \Psi_{\varepsilon,n},$$

where $|\Psi_{\varepsilon,n}| \le c(n)\varepsilon$.

The process $\eta$ has only a finite number of states, so $L\hat{u}_n$ is uniformly bounded. The process $\hat{u}_n$ jumps after an exponential waiting time. These facts and the formula $M_n(t) = \hat{u}_n(t) - \int_0^t L\hat{u}_n(s)\,ds$ imply that

$$\lim_{s\to 0} \frac{E[(M_n(t+s) - M_n(t))^2]}{E[(\hat{u}_n(t+s) - \hat{u}_n(t))^2]} = 1.$$

It follows that

$$E[M_n^2(t)] = \int_0^t E[A_n(s)]\,ds, \tag{3.8}$$

where $A_n = A_n(t, \eta)$ is given by the following formula:

$$A_n = \lim_{s\to 0} (1/s) E[(\hat{u}_n(t+s) - \hat{u}_n(t))^2 | \eta(t)]$$

$$= N^{-2} \sum_{x,y\in D_\varepsilon} h_\varepsilon^{-1}(x) p_{xy} \Bigg\{ \eta_x^+ \mathbf{1}_{\{\eta_y \ge 0\}} (\phi_n(y) - \phi_n(x))^2$$

$$+ \eta_x^- \mathbf{1}_{\{\eta_y \le 0\}} (\phi_n(x) - \phi_n(y))^2$$

$$+ \eta_x^+ \mathbf{1}_{\{\eta_y < 0\}} N^{-2} \sum_{u,v\in D_\varepsilon} \eta_u^+ \eta_v^- (\phi_n(y) - \phi_n(x)$$

$$+ \phi_n(u) - \phi_n(v))^2$$

$$+ \eta_x^- \mathbf{1}_{\{\eta_y > 0\}} N^{-2} \sum_{u,v\in D_\varepsilon} \eta_u^- \eta_v^+ (\phi_n(x) - \phi_n(y)$$

$$- \phi_n(u) + \phi_n(v))^2 \Big\}.$$

This implies the following bound for $A_n$, for small $\varepsilon$:

$$A_n \le c_1 N^{-1} V(\varepsilon^2 \|\nabla \phi_n\|_\infty^2 + \|\phi_n\|_\infty^2) \le \beta_n N^{-1} V, \tag{3.9}$$

where $\beta_n < \infty$ depends on $\phi_n$, and $\nabla$ stands for the usual gradient acting on functions defined on $\mathbb{R}^d$. It follows from (3.6) and (3.7) that

$$\begin{aligned}
\hat{u}_n(t) &= e^{\int_0^t (V(r)-\lambda_n)\,dr} \Big( \hat{u}_n(0) + \int_0^t e^{-\int_0^s (V(r)-\lambda_n)dr}\, dM_s \\
&\qquad + \int_0^t e^{-\int_0^s (V(r)-\lambda_n)dr} \Psi_{\varepsilon,n}(s)\, ds \Big) \\
&= \hat{u}_n(0) e^{\int_0^t (V(r)-\lambda_n)\,dr} + \int_0^t e^{\int_s^t (V(r)-\lambda_n)\,dr}\, dM_s \\
&\quad + \int_0^t e^{\int_s^t (V(r)-\lambda_n)\,dr} \Psi_{\varepsilon,n}(s)\, ds \\
&\stackrel{\mathrm{df}}{=} \hat{u}_n(0) e^{\int_0^t (V(r)-\lambda_n)\,dr} + R_{n,1}(t) + R_{n,2}(t).
\end{aligned} \tag{3.10}$$

Since $R_{n,1}(t) = \int_0^t e^{\int_s^t (V(r)-\lambda_n)\,dr}\, dM_s$, (3.8) and [20], Chapter II, Section 6, Corollary 3 and Theorem 29, show that

$$E[R_{n,1}^2(t)] = E\Big[\int_0^t e^{2\int_s^t (V(r)-\lambda_n)\,dr} A_n(s)\, ds\Big].$$

In view of (3.9),

$$E[R_{n,1}^2(t)] \le \beta_n N^{-1} E\Big[\int_0^t V(s) e^{2\int_s^t (V(r)-\lambda_n)\,dr}\, ds\Big].$$

We have

$$\begin{aligned}
&\int_0^t V(s) e^{2\int_s^t (V(r)-\lambda_n)\,dr}\, ds \\
&\quad = \lambda_n \int_0^t e^{2\int_s^t (V(r)-\lambda_n)\,dr}\, ds + (1/2)(e^{2\int_0^t (V(r)-\lambda_n)\,dr} - 1) \\
&\quad \le c_2 e^{2\int_0^t V(r)\,dr},
\end{aligned}$$

where $c_2$ depends on $t$ and $n$. Thus, for some $c_3$ that depends on $t$ and $n$,

$$E[R_{n,1}^2(t)] \le c_3 \beta_n N^{-1} E\Big[\int_0^t e^{2\int_s^t V(r)\,dr}\, ds\Big].$$

For the second remainder, we have the following estimate, using the Cauchy–Schwarz inequality:

$$\begin{aligned} E[R_{n,2}^2(t)] &\le t^2\Big(\sup_{s\le t}\Psi^2_{\varepsilon,n}(s)\Big)E\bigg[\int_0^t e^{2\int_s^t (V(r)-\lambda_n)\,dr}\,ds\bigg] \\ &\le c(n)^2\varepsilon^2 t^2 E\bigg[\int_0^t e^{2\int_s^t (V(r)-\lambda_n)\,dr}\,ds\bigg]. \end{aligned}$$

Hence,

$$\begin{aligned} &E[(R_{n,1}(t)+R_{n,2}(t))^2] \\ &\qquad \le 2(c_3\beta_n N^{-1} + c(n)^2\varepsilon^2 t^2)E[e^{2\int_0^t V(r)\,dr}]. \end{aligned} \tag{3.11}$$

Recall that, for some $n$, there exists $a>0$ such that $\inf_{N,\varepsilon}|\hat{u}_n(0)| = a$. Let $n_0$ be the smallest $n$ satisfying this condition. It follows from (3.10) and (3.11) that

$$\begin{aligned} &E[(\hat{u}_{n_0}(t) - \hat{u}_{n_0}(0)e^{\int_0^t (V(r)-\lambda_{n_0})\,dr})^2] \\ &\qquad \le 2(c_3\beta_{n_0} N^{-1} + c(n_0)^2\varepsilon^2 t^2)E[e^{2\int_0^t V(r)\,dr}]. \end{aligned} \tag{3.12}$$

Suppose that $N$ is large enough (and, consequently, $\varepsilon$ is small) so that

$$4(c_3\beta_{n_0}N^{-1} + c(n_0)^2\varepsilon^2 t^2) \le (1/2)a^2 e^{-2\lambda_{n_0} t}. \tag{3.13}$$

Since $\|\phi_{n_0}\|_\infty < \infty$ and $\varepsilon^d \sum_x |u(x,t)| = 2$, we obtain

$$|\hat{u}_{n_0}(t)| = \bigg|\varepsilon^d \sum_{x\in D_\varepsilon} u(x,t)\phi_{n_0}(x)\bigg| \le c_4. \tag{3.14}$$

By (3.12), (3.13) and (3.14),

$$\begin{aligned} &a^2 e^{-2\lambda_{n_0} t} E[e^{2\int_0^t V(r)\,dr}] \\ &\quad \le E[(\hat{u}_{n_0}(0)e^{\int_0^t (V(r)-\lambda_{n_0})\,dr})^2] \\ &\quad \le 2E[(\hat{u}_{n_0}(t) - \hat{u}_{n_0}(0)e^{\int_0^t (V(r)-\lambda_{n_0})\,dr})^2] + 2E[(\hat{u}_{n_0}(t))^2] \\ &\quad \le 4(c_3\beta_{n_0}N^{-1} + c(n_0)^2\varepsilon^2 t^2)E[e^{2\int_0^t V(r)\,dr}] + 2c_4^2 \\ &\quad \le (1/2)a^2 e^{-2\lambda_{n_0} t}E[e^{2\int_0^t V(r)\,dr}] + 2c_4^2, \end{aligned}$$

so

$$E[e^{2\int_0^t V(r)\,dr}] \le 4c_4^2 a^{-2} e^{2\lambda_{n_0} t}. \tag{3.15}$$

Note that $\widehat{P_t u(0)}_n = e^{-\lambda_n t}\hat{u}_n(0)$ and let $v(x,t) = e^{\int_0^t V(r)dr}(P_t u^{N,\varepsilon}(0))(x)$. We combine (3.10), (3.11) and (3.15) to see that if $N \to \infty$ (and, therefore, $\varepsilon \to 0$), then, for every fixed $n$ and $t > 0$, $\hat{u}_n(t) - \hat{v}_n(t) \to 0$ in distribution. Then Lemma 2.2 shows that

$$\sum_{x\in D_\varepsilon} \varepsilon^d u^{N,\varepsilon}(x,t)\mathbf{i}_x - e^{\int_0^t V(r)\,dr} P_t u^{N,\varepsilon}(0) \to \delta_{\mathbf{0}} \qquad \text{in } M_1(M_F(\overline{D})). \tag{3.16}$$

*Step* 3. It is not obvious that the normalization of $P_t u^{N,\varepsilon}(0)$ in (3.16) is the same as in the statement of Theorem 3.2. It is conceivable that a sizeable proportion of positive and negative particles are tightly interspersed so that their masses cancel each other in the limit. We will show that this is not the case—intuitively, the two populations occupy disjoint parts of $D$.

Let $B_\delta(x)$ denote a hypercube in $D_\varepsilon$, centered at $x$, with side length $\delta$. We will consider only $\delta > \varepsilon$. We set

$$\Lambda_\delta(t,x) = \min\left(\sum_{y\in B_\delta(x)} \eta_y^+, \sum_{y\in B_\delta(x)} \eta_y^-\right).$$

Note that if $\varepsilon < \delta/2$ and $x \in D_\delta$, then $B_\delta(x)$ contains at least $(\delta/2\varepsilon)^d$ sites $y \in D_\varepsilon$. Fix an arbitrarily small $c_0 > 0$. Suppose that, for some $B_\delta(x)$, we have $\Lambda_\delta(t,x) \ge c_0(\delta/2\varepsilon)^d$. Given this assumption, we will show that the number of ($+$ and $-$) particles that are located in $B_\delta(x)$ at time $t$ and collide with a particle of the opposite sign before time $t+\delta^2$ has expectation greater than $c_1\Lambda_\delta(t,x)$. Suppose, without loss of generality, that there are fewer $+$ than $-$ particles in $B_\delta(x)$ at time $t$, that is, $\Lambda_\delta(t,x) = \sum_{y\in B_\delta(x)} \eta_x^+$. Consider independent continuous time reflecting random walks $Y_k$, $1 \le k \le \Lambda_\delta(t,x)$, starting from the same points as the locations of the $+$ particles at time $t$ in $B_\delta(x)$. The distribution of a process $\{Y_k(s), s \ge t\}$ is the same as for a single particle in our process $\eta$, except that $Y_k$'s do not interact with other particles. For $k < (1/2)(\delta/2\varepsilon)^d \wedge \Lambda_\delta(t,x)$,

$$P(Y_k(t+\delta^2) \ne Y_j(t+\delta^2), 1 \le j \le k-1) > p_1 > 0,$$

where $p_1$ depends only on the dimension $d$. This implies that, with probability $p_2$, the number of different sites occupied by $Y_k(t+\delta^2)$, $1 \le k \le \Lambda_\delta(t,x)$, is greater than $c_2(\delta/\varepsilon)^d$, where $p_2, c_2 > 0$ depend only on $d$ and $c_0$.

Choose $\Lambda_\delta(t,x)$ locations of the $-$ particles at time $t$ in $B_\delta(x)$. From each of these points, start a continuous time reflecting random walk $Z_k$. Assume that $Z_k$'s, $1 \le k \le \Lambda_\delta(t,x)$, are independent, and they are independent of $Y_k$'s. Suppose that the number of different sites occupied by $Y_k(t+\delta^2)$, $1 \le k \le \Lambda_\delta(t,x)$, is greater than $c_2(\delta/\varepsilon)^d$ and call the set of these sites $\Gamma$. Then it is easy to see that, with probability greater than $p_3 > 0$, the number of distinct sites in $\Gamma$ occupied by $Z_k(t+\delta^2)$, $1 \le k \le (c_2/2)(\delta/\varepsilon)^d$, is greater than $c_3(\delta/\varepsilon)^d$.

Suppose that there are at least $c_2(\delta/\varepsilon)^d$ sites in $\Gamma$ and the number of distinct sites in $\Gamma$ occupied by $Z_k(t+\delta^2)$, $1 \le k \le (c_2/2)(\delta/\varepsilon)^d$, is greater than $c_3(\delta/\varepsilon)^d$. Find the first time $t_1 > t$ when some $Y_k$ and $Z_j$ occupy the same site and call these particles "eliminated." Then, by induction, find the smallest $t_m > t_{m-1}$ when some noneliminated $Y_k$ and $Z_j$ occupy the same site and eliminate this pair of particles. Note that the total number of eliminated pairs by the time $t+\delta^2$ cannot be smaller than $(c_3/2)(\delta/\varepsilon)^d$.

Now we return to our original model, with interactions between particles. Consider the set of $+$ and $-$ particles in $\eta$ that reside at the same locations as $Y_k$'s and $Z_k$'s at time $t$ in $B_\delta(x)$. Choose from this set a pair $(Q_+, Q_-)$ consisting of a $+$ and a $-$ particle and suppose that it would have been "eliminated" in the scheme described above, that is, if these two particles had been $Y_k$ and $Z_j$ for some $k$ and $j$. If $Q_+$ and $Q_-$ do not meet before time $t+\delta^2$, it means that one of these particles must have met a particle of the opposite sign (different from $Q_+$ and $Q_-$) before time $t+\delta^2$, and hence, at least one of particles $Q_+$ and $Q_-$ has a jump before time $t+\delta^2$. Thus, with probability greater than $p_2p_3 > 0$, there will be at least $(c_3/4)(\delta/\varepsilon)^d$ jumps between times $t$ and $t+\delta^2$, by particles that are located in $B_\delta(x)$ at time $t$, assuming that $\Lambda_\delta(t,x) \ge c_0(\delta/2\varepsilon)^d$.

Let $K(t)$ be the number of collisions before $t$. Then $EK(t) = (N/2)E\int_0^t V(s)\,ds$. Let $H_\delta(t) = \{x \in D_\delta : \Lambda_\delta(t,x) \ge c_0(\delta/2\varepsilon)^d\}$. We see that, for some $c_1$ depending on $c_0$,

$$EK(t) \ge \sum_{k=0}^{[t/\delta^2]-1} E \sum_{x\in H_\delta(k\delta^2)} c_1[\Lambda_\delta(k\delta^2, x)].$$

For the same reason, for $s \in [0, \delta^2]$,

$$EK(t) \ge \sum_{k=0}^{[(t-s)/\delta^2]-1} E \sum_{x\in H_\delta(k\delta^2)} c_1[\Lambda_\delta(k\delta^2 + s, x)].$$

Hence,

$$EK(t) \ge \int_{\delta^2}^{t-2\delta^2} \delta^{-2} E \sum_{x\in H_\delta(s)} c_1[\Lambda_\delta(s,x)]\,ds. \tag{3.17}$$

For a function $f$ on $D_\varepsilon$, let $\Xi(x,\delta,f) = \sum_{y\in B_\delta(x)} f(y)$. Recall that $u_t = u(x,t) = N^{-1}\varepsilon^{-d}\eta_x$. We have

$$\Xi(x,\delta,|u_t|) - |\Xi(x,\delta,u_t)| = 2N^{-1}\varepsilon^{-d}\Lambda_\delta(t,x). \tag{3.18}$$

If $x \notin H_\delta(s)$, then either

$$\sum_{y\in B_\delta(x)} \eta_y^+ \le c_0(\delta/2\varepsilon)^d \quad \text{or} \quad \sum_{y\in B_\delta(x)} \eta_y^- \le c_0(\delta/2\varepsilon)^d.$$

Hence,

$$\tag{3.19}\begin{aligned}\sum_{x\notin H_\delta(s)}(\Xi(x,\delta,|u_s|)-|\Xi(x,\delta,u_s)|)&\le\sum_{x\notin H_\delta(s)}c_0(\delta/2\varepsilon)^d\\&\le\sum_{x\in D_\delta}c_0(\delta/2\varepsilon)^d\le c_0c_4\varepsilon^{-d},\end{aligned}$$

where $c_4$ is a constant depending only on $D$. Recall that $N\ge c_5\varepsilon^{-d}$. In view of (3.15), (3.17), (3.18) and (3.19),

$$\tag{3.20}\begin{aligned}&E\left[\int_{\delta^2}^{t-2\delta^2}\left(2-\sum_{x\in D_\delta}|\Xi(x,\delta,\varepsilon^d u_s)|\right)ds\right]\\&\quad=E\left[\int_{\delta^2}^{t-2\delta^2}\sum_{x\in D_\delta}(\Xi(x,\delta,|\varepsilon^d u_s|)-|\Xi(x,\delta,\varepsilon^d u_s)|)\,ds\right]\\&\quad=\varepsilon^d E\left[\int_{\delta^2}^{t-2\delta^2}\sum_{x\in D_\delta}(\Xi(x,\delta,|u_s|)-|\Xi(x,\delta,u_s)|)\,ds\right]\\&\quad=\varepsilon^d E\left[\int_{\delta^2}^{t-2\delta^2}\sum_{x\notin H_\delta(s)}(\Xi(x,\delta,|u_s|)-|\Xi(x,\delta,u_s)|)\,ds\right]\\&\qquad+\varepsilon^d E\left[\int_{\delta^2}^{t-2\delta^2}\sum_{x\in H_\delta(s)}(\Xi(x,\delta,|u_s|)-|\Xi(x,\delta,u_s)|)\,ds\right]\\&\quad\le c_0c_4t+2\varepsilon^d E\left[\int_{\delta^2}^{t-2\delta^2}\sum_{x\in H_\delta(s)}N^{-1}\varepsilon^{-d}\Lambda_\delta(s,x)\,ds\right]\\&\quad\le c_0c_4t+2c_6N^{-1}\delta^2EK(t)\\&\quad=c_0c_4t+2c_6\delta^2E\left[\int_0^t V(s)\,ds\right]\\&\quad\le c_0c_4t+c_7\delta^2,\end{aligned}$$

where $c_7$ depends on $c_0$. For a fixed $t>0$ and $n>0$, we can find $c_0>0$ and $\delta_n\in(0,2^{-n})$ so small that the right-hand side of (3.20) is less than $2^{-n}$. Let $\mathcal{T}_n$ be the set of $s\in[0,t]$ such that, for $\varepsilon<\delta_n$,

$$E\left[2-\sum_{x\in D_{\delta_n}}|\Xi(x,\delta,\varepsilon^d u_s)|\right]\ge n^2 2^{-n}.$$

Then, by (3.20), $|\mathcal{T}_n| \le n^{-2}$. Choose arbitrarily small $c_* > 0$ and let $n_0$ be so large that $\sum_{n\ge n_0} |\mathcal{T}_n| < c_*$. Let $\mathcal{T}_* = (0,t] \setminus \bigcup_{n\ge n_0} \mathcal{T}_n$. For $s \in \mathcal{T}_*$, and $\varepsilon < \delta_n$,

$$P\left(2 - \sum_{x\in D_{\delta_n}} |\Xi(x,\delta,\varepsilon^d u_s)| \ge n^3 2^{-n}\right) \le n^{-1}. \tag{3.21}$$

By passing to a subsequence, if necessary, we may assume that $\sum_{x\in D_\varepsilon} \varepsilon^d u^{N,\varepsilon}(x, s)\mathbf{i}_x$ converges in $M_1(M_F(\overline{D}))$. It follows from (3.21) that any limit of $\sum_{x\in D_\varepsilon} \varepsilon^d u^{N,\varepsilon}(x,s)\mathbf{i}_x$ is supported on measures with the total variation 2, for every $s \in \mathcal{T}_*$. Since $c_* > 0$ in the definition of $\mathcal{T}_*$ is arbitrarily small, we obtain the same conclusion for almost every $s \in [0,t]$. This and (3.16) imply that, for almost every $s > 0$,

$$\sum_{x\in D_\varepsilon} \varepsilon^d u^{N,\varepsilon}(x,s)\mathbf{i}_x - \overline{P_s u^{N,\varepsilon}(0)} \to \delta_{\mathbf{0}} \qquad \text{in } M_1(M_F(\overline{D})). \tag{3.22}$$

*Step* 4. We will next show that, for large $N$, the process $r \to \int_0^r V(s)\,ds$ is close to being continuous, in the sense that

$$\limsup_{\delta_1,\delta_2\to 0} \limsup_{\varepsilon\to 0} P\left(\sup_{0\le t_1,t_2\le t, |t_1-t_2|\le\delta_1} \int_{t_1}^{t_2} V(s)\,ds \ge \delta_2\right) = 0.$$

We need the above claim in the last step of the proof.

Let $C(s)$ be such that $C(s)P_s u^{N,\varepsilon}(0) = \overline{P_s u^{N,\varepsilon}(0)}$. Note that $C(s) \in (0,\infty)$ for all $s$ because the solution to the heat equation in $\overline{D}$ with Neumann boundary conditions and nonzero initial condition with finite total variation is never identically zero and for all $s$, it has a finite variation. It is easy to see that $C(s)$ is a nondecreasing function. By assumption, $C(0) = 1$. Recall that, for some $n$, there exists $a > 0$ such that $\inf_{N,\varepsilon} |\hat{u}_n(0)| = a$ and $\widehat{P_s u(0)}_n = e^{-\lambda_n s}\hat{u}_n(0)$. Since $\widehat{P_s u(0)}_n = \int \phi_n\, dP_s u(0)$ and $\phi_n$ is bounded (see the proof of Lemma 2.2), it follows that the total variation of $P_t u(0)$ is bounded below and, therefore, $C(t)$ is bounded above by a constant depending on $D, n, a$ and $t$.

By Theorem 2.1 of [3], for any $0 < s_1 < s_2 < \infty$, the Neumann heat kernel $p_s(x,y)$ is Hölder continuous jointly in $(s,x,y)$ on $[s_1,s_2] \times \overline{D} \times \overline{D}$. Let $f(s,x)$ be the density of $P_s u(0)$ at $x \in \overline{D}$. We see that, for any $0 < s_1 < s_2 < \infty$, the family of densities $\{f(s,x)\}$ corresponding to all measures $u^{N,,\varepsilon}(0)$ satisfying our assumptions is equicontinuous on $[s_1,s_2] \times \overline{D}$. This implies that the family of functions $\{C(s), 0 \le s \le t\}$ is equicontinuous. Since the family $\{f(s,x)\}$ is equicontinuous and the total variation of $P_s u(0)$ is bounded below on $[s_1,s_2] \times \overline{D}$, it follows that the Prohorov distance of $P_s u(0)$ from the measure identically equal to 0 is bounded below on $[s_1,s_2]$ by a constant depending only on $D, n, a, s_1$ and $s_2$. By (3.16) and (3.22), for almost every $s > 0$,

$$(C(s) - e^{\int_0^s V(r)\,dr})P_s u^{N,\varepsilon}(0) \to \delta_{\mathbf{0}} \qquad \text{in } M_1(M_F(\overline{D})).$$

We see that $C(s) - e^{\int_0^s V(r)\,dr}$ must converge to 0 in distribution for almost every $s$. Since $C(s)$ are uniformly continuous and nondecreasing, and $e^{\int_0^s V(r)\,dr}$ is nondecreasing, it is easy to see that $C(s) - e^{\int_0^s V(r)\,dr}$ converges to 0 uniformly on compact intervals, in distribution, as $N \to \infty$. This implies the claim stated at the beginning of Step 4.

*Step* 5. Next we will show that the convergence holds not only for every fixed $t > 0$, but also in the Skorohod topology on $M_1(S((0,\infty), M_F(\overline{D})))$. Fix a smooth function $\varphi$ on $\overline{D}$ and let $w_t^\varepsilon = w_t = \langle u_t, \varphi\rangle$. We will show that the family of processes $\{w^\varepsilon, \varepsilon > 0\}$ is tight in $M_1(S((0,\infty), \mathbb{R}))$. In order to prove that, we will first derive some estimates for $w_t$ similar to the estimates for $\hat{u}_n$. We have

$$dw_t = Lw_t\,dt + dM_\varphi(t), \tag{3.23}$$

where $M_\varphi(t)$ is a martingale. Since $\varphi$ is smooth, there exists a series expansion for $\varphi$ that yields

$$\Delta_\varepsilon \varphi(x) = \Delta\varphi + \psi_\varepsilon(x),$$

where $|\psi_\varepsilon(x)| \le c_\varphi \varepsilon$. We have $L\eta_x = \Delta_\varepsilon^* \eta_x + V\eta_x$ and

$$w_t = \langle u, \varphi\rangle = \varepsilon^d \sum_{x\in D_\varepsilon} N^{-1}\varepsilon^{-d}\eta_x\varphi(x) = N^{-1}\sum_{x\in D_\varepsilon}\eta_x\varphi(x),$$

so

$$\begin{aligned}
Lw &= N^{-1}\sum_{x\in D_\varepsilon}(L\eta_x)\varphi(x) \\
&= N^{-1}\sum_{x\in D_\varepsilon}(\Delta_\varepsilon^*\eta_x)\varphi(x) + N^{-1}\sum_{x\in D_\varepsilon}(V\eta_x)\varphi(x) \\
&= N^{-1}\sum_{x\in D_\varepsilon}\left(\sum_{y\in D_\varepsilon}(h_\varepsilon^{-1}(y)p_{yx}\eta_y - h_\varepsilon^{-1}(x)p_{xy}\eta_x)\right)\varphi(x) \\
&\quad + N^{-1}V\sum_{x\in D_\varepsilon}\eta_x\varphi(x) \\
&= N^{-1}\sum_{x\in D_\varepsilon}\left(\sum_{y\in D_\varepsilon}(h_\varepsilon^{-1}(x)p_{xy}\varphi(y) - h_\varepsilon^{-1}(x)p_{xy}\varphi(x))\right)\eta_x \\
&\quad + N^{-1}V\sum_{x\in D_\varepsilon}\eta_x\varphi(x)
\end{aligned}$$

$$= N^{-1}\sum_{x\in D_\varepsilon}\eta_x(\Delta_\varepsilon\varphi(x)) + VN^{-1}\sum_{x\in D_\varepsilon}\eta_x\varphi(x) \tag{3.24}$$

$$= N^{-1} \sum_{x \in D_\varepsilon} \eta_x (\Delta\varphi(x) + \psi_\varepsilon(x)) + V N^{-1} \sum_{x \in D_\varepsilon} \eta_x \varphi(x)$$
$$= \langle u, \Delta\varphi \rangle + \Psi_\varepsilon + V w,$$

where $|\Psi_\varepsilon| \le c_\varphi \varepsilon$. Note that

$$E[M_\varphi^2(t)] = \int_0^t E[A_\varphi(s)]\, ds, \tag{3.25}$$

where $A_\varphi = A_\varphi(t, \eta)$ is given by the following formula:

$$A_\varphi = \lim_{s \to 0} (1/s) E[(w(t+s) - w(t))^2 | \eta(t)]$$
$$= N^{-2} \sum_{x,y \in D_\varepsilon} h_\varepsilon^{-1}(x) p_{xy} \Big\{ \eta_x^+ \mathbf{1}_{\{\eta_y \ge 0\}} (\varphi(y) - \varphi(x))^2$$
$$+ \eta_x^- \mathbf{1}_{\{\eta_y \le 0\}} (\varphi(x) - \varphi(y))^2$$
$$+ \eta_x^+ \mathbf{1}_{\{\eta_y < 0\}} N^{-2} \sum_{u,v \in D_\varepsilon} \eta_u^+ \eta_v^- (\varphi(y) - \varphi(x)$$
$$+ \varphi(u) - \varphi(v))^2$$
$$+ \eta_x^- \mathbf{1}_{\{\eta_y > 0\}} N^{-2} \sum_{u,v \in D_\varepsilon} \eta_u^- \eta_v^+ (\varphi(x) - \varphi(y)$$
$$- \varphi(u) + \varphi(v))^2 \Big\}.$$

We obtain the following bound for $A_\varphi$, for small $\varepsilon$:

$$A_\varphi \le c_1 N^{-1} V (\varepsilon^2 \|\nabla\varphi\|_\infty^2 + \|\varphi\|_\infty^2) \le \beta N^{-1} V, \tag{3.26}$$

where $\beta < \infty$ depends on $\varphi$.

It follows from (3.23) and (3.24) that

$$w(t) = e^{\int_0^t V(r)\, dr} \Big( w(0) + \int_0^t e^{-\int_0^s V(r)\, dr}\, dM_\varphi(s)$$
$$+ \int_0^t e^{-\int_0^s V(r)\, dr} (\langle u, \Delta\varphi \rangle(s) + \Psi_\varepsilon(s))\, ds \Big)$$
$$= w(0) e^{\int_0^t V(r)\, dr} + \int_0^t e^{\int_s^t V(r)\, dr}\, dM_\varphi(s) \tag{3.27}$$
$$+ \int_0^t e^{\int_s^t V(r)\, dr} \Psi_\varepsilon(s)\, ds + \int_0^t e^{\int_s^t V(r)\, dr} \langle u, \Delta\varphi \rangle(s)\, ds$$
$$\stackrel{\text{df}}{=} w(0) e^{\int_0^t V(r)\, dr} + R_1(t) + R_2(t) + R_3(t).$$

In view of (3.25), (3.26) and (3.15), we have

$$
\begin{aligned}
E[R_1^2(t)] &= E\left[\int_0^t e^{2\int_s^t V(r)\,dr} A_n(s)\,ds\right] \\
&\le c_2\beta N^{-1} E\left[\int_0^t V(s) e^{2\int_s^t V(r)\,dr}\,ds\right] \\
&= c_2\beta N^{-1} E[(1/2)(e^{2\int_0^t V(r)\,dr} - 1)] \\
&\le c_3\beta N^{-1} E[e^{2\int_0^t V(r)\,dr}] \\
&\le c_4 e^{2\lambda_{n_0} t} N^{-1}.
\end{aligned}
$$

Since $R_1(t) = \int_0^t e^{\int_s^t V(r)\,dr}\,dM_\varphi(s)$ and $M_\varphi(t)$ is a martingale, so is $R_1(t)$. Hence, by Doob's inequality,

$$
E\left[\sup_{0\le s\le t} R_1^2(s)\right] \le 4E[R_1^2(t)] \le 4c_4 e^{2\lambda_{n_0} t} N^{-1}.
$$

A similar calculation and (3.15) show that

$$
E\left[\sup_{0\le s\le t} R_2^2(s)\right] \le t^2\left(\sup_{0\le s\le t} \Psi_\varepsilon^2(s)\right) E\left[\int_0^t e^{2\int_s^t V(r)\,dr}\,ds\right] \le c_5\varepsilon^2 t^2.
$$

Hence,

$$
E\left[\sup_{0\le s\le t} (R_1(s) + R_2(s))^2\right] \le c_6(N^{-1} + \varepsilon^2 t^2), \tag{3.28}
$$

where $c_6$ depends on $t$ and $\varphi$.

For the last term on the right-hand side of (3.27), we observe that, for $0 \le t_1 < t_2 \le t$,

$$
\begin{aligned}
|R_3(t_2) - R_3(t_1)| &= \bigg| \int_{t_1}^{t_2} e^{\int_s^{t_2} V(r)\,dr} \langle u, \Delta\varphi\rangle(s)\,ds \\
&\qquad + \int_0^{t_1} (e^{\int_s^{t_2} V(r)\,dr} - e^{\int_s^{t_1} V(r)\,dr}) \langle u, \Delta\varphi\rangle(s)\,ds \bigg| \\
&= \bigg| \int_{t_1}^{t_2} e^{\int_s^{t_2} V(r)\,dr} \langle u, \Delta\varphi\rangle(s)\,ds \\
&\qquad + \int_0^{t_1} e^{\int_s^{t_1} V(r)\,dr} (e^{\int_{t_1}^{t_2} V(r)\,dr} - 1) \langle u, \Delta\varphi\rangle(s)\,ds \bigg| \\
&\le \sup_{0\le s\le t} e^{\int_s^t V(r)\,dr} \int_{t_1}^{t_2} |\langle u, \Delta\varphi\rangle(s)|\,ds \\
&\qquad + (e^{\int_{t_1}^{t_2} V(r)\,dr} - 1) \sup_{0\le s\le t} e^{\int_s^t V(r)\,dr} \int_0^{t_1} |\langle u, \Delta\varphi\rangle(s)|\,ds
\end{aligned}
$$

$$\le c_7 e^{\int_0^t V(r)\,dr} \|\Delta\varphi\|_\infty |t_2 - t_1|$$
$$+ c_8 (e^{\int_{t_1}^{t_2} V(r)\,dr} - 1) t e^{\int_0^t V(r)\,dr} \|\Delta\varphi\|_\infty.$$

In view of (3.15) and the first claim in Step 4,

$$\limsup_{\delta_1,\delta_2\to 0} \limsup_{\varepsilon\to 0} P\left( \sup_{0\le t_1,t_2\le t, |t_1-t_2|\le\delta_1} |R_3(t_2) - R_3(t_1)| \ge \delta_2 \right) = 0.$$

This, (3.27), (3.28) and Lemma 2.3 show that $\{w^\varepsilon, \varepsilon > 0\}$ is a tight family of processes in $M_1(S((0,\infty),\mathbb{R}))$. Since smooth functions are dense in the set of continuous functions on $\overline{D}$ and the sum of two smooth functions is smooth, Theorem 3.7.1 of [12] shows that the family of processes $\{\sum_{x\in D_\varepsilon} \varepsilon^d u^{N,\varepsilon}(x,\cdot)\mathbf{i}_x, \varepsilon > 0\}$ is tight in $M_1(S((0,\infty), M_F(\overline{D})))$. Since $\{P_0 u^{N,\varepsilon}(0), \varepsilon > 0\}$ are tight and the process $t \to \overline{P_t u^{N,\varepsilon}(0)}$ is a continuous function with values in $M_F(\overline{D})$, completely determined by $P_0 u^{N,\varepsilon}(0)$, we conclude that

$$\left\{ \sum_{x\in D_\varepsilon} \varepsilon^d u^{N,\varepsilon}(x,\cdot)\mathbf{i}_x - \overline{P_\cdot u^{N,\varepsilon}(0)}, \varepsilon > 0 \right\}$$

is tight in $M_1(S((0,\infty), M_F(\overline{D})))$. By (2.13), any convergent subsequence of this family must be identically equal to $\delta_{\mathbf{0}}$, by the right continuity. □

**Acknowledgments.** We are grateful to Rodrigo Bañuelos, Victor Ivrii, John Sylvester, Alain Sznitman, John Toth and Steve Zelditch for very helpful advice. We thank the referees for many helpful suggestions.

## REFERENCES

[1] Bañuelos, R. and Burdzy, K. (1999). On the "hot spots" conjecture of J. Rauch. *J. Func. Anal.* **164** 1–33. MR1694534

[2] Bañuelos, R. and Pang, M. (2006). Level sets of Neumann eigenfunctions. *Indiana Univ. Math. J.* **55** 923–939. MR2244591

[3] Burdzy, K. and Chen, Z. (2002). Coalescence of synchronous couplings. *Probab. Theory Related Fields* **123** 553–578. MR1921013

[4] Burdzy, K., Hołyst, R., Ingerman, D. and March, P. (1996). Configurational transition in a Fleming–Viot-type model and probabilistic interpretation of Laplacian eigenfunctions. *J. Phys. A* **29** 2633–2642.

[5] Burdzy, K., Hołyst, R. and March, P. (2000). A Fleming–Viot particle representation of Dirichlet Laplacian. *Comm. Math. Phys.* **214** 679–703. MR1800866

[6] Conti, M., Terracini, S. and Verzini, G. (2005). Asymptotic estimates for the spatial segregation of competitive systems. *Adv. Math.* **195** 524–560. MR2146353

[7] Conti, M., Terracini, S. and Verzini, G. (2005). A variational problem for the spatial segregation of reaction–diffusion systems. *Indiana Univ. Math. J.* **54** 779–815. MR2151234

[8] Courant, R. and Hilbert, D. (1953). *Methods of Mathematical Physics*. Interscience Publishers, New York. MR0065391

[9] Cybulski, O., Babin, V. and Hołyst, R. (2004). Minimization of the Renyi entropy production in the stationary states of the Brownian process with matched death and birth rates. *Phys. Rev. E* **69** 016110.

[10] Cybulski, O., Babin, V. and Hołyst, R. (2005). Minimization of the Renyi entropy production in the space partitioning process. *Phys. Rev. E* **71** 046130. MR2139992

[11] Cybulski, O., Matysiak, D., Babin, V. and Hołyst, R. (2005). Pattern formation in nonextensive thermodynamics: Selection criterion based on the Renyi entropy production. *J. Chem. Phys.* **122** 174105.

[12] Dawson, D. A. (1993). Measure-valued Markov processes. *École d'Été de Probabilités de Saint-Flour XXI—1991* 1–260. *Lecture Notes in Math.* **1541**. Springer, Berlin. MR1242575

[13] Ethier, S. N. and Kurtz, T. G. (1986). *Markov Processes. Characterization and Convergence*. Wiley, New York. MR0838085

[14] Grieser, D. (2002). Uniform bounds for eigenfunctions of the Laplacian on manifolds with boundary *Comm. Partial Diff. Eq.* **27** 1283–1299. MR1924468

[15] Grigorescu, I. and Kang, M. (2004). Hydrodynamic limit for a Fleming–Viot type system. *Stochastic Process. Appl.* **110** 111–143. MR2052139

[16] Grigorescu, I. and Kang, M. (2006). Tagged particle limit for a Fleming–Viot type system. *Electron. J. Probab.* **11** 311–331. MR2217819

[17] Hempel, R., Seco, L. A. and Simon, B. (1991). The essential spectrum of Neumann Laplacians on some bounded singular domains. *J. Func. Anal.* **102** 448–483. MR1140635

[18] Morrey, C. (1966). *Multiple Integrals in the Calculus of Variations*. Springer, New York. MR0202511

[19] Netrusov, Yu. and Safarov, Yu. (2005). Weyl asymptotic formula for the Laplacian on domains with rough boundaries. *Comm. Math. Phys.* **253** 481–509. MR2140257

[20] Protter, Ph. (1990). *Stochastic Integration and Differential Equations. A New Approach*. Springer, Berlin. MR1037262

[21] Stroock, D. and Varadhan, S. R. S. (1971). Diffusion processes with boundary conditions. *Comm. Pure Appl. Math.* **24** 147–225. MR0277037

Department of Mathematics
University of Washington
Box 354350
Seattle, Washington 98115-4350
USA
E-mail: burdzy@math.washington.edu

Department of Mathematics
University of Toronto
100 St. George St.
Toronto, Ontario
Canada M5S 3G3
E-mail: quastel@math.toronto.edu